\documentstyle[amstex,amscd,pb-diagram,leqno]{article}

\newtheorem{thm}{Theorem}[section]

\newtheorem{lem}[thm]{Lemma}
\newtheorem{cor}[thm]{Corollary}
\newtheorem{defn}[thm]{Definition}
\newtheorem{prop}[thm]{Proposition}
\newtheorem{remk}[thm]{Remark}
\newtheorem{noten}[thm]{Notation}
\newtheorem{exam}[thm]{Example}

\newtheorem{prob}[thm]{Problem}
\newtheorem{whypo}[thm]{Working Hypothesis}
\newtheorem{conj}[thm]{Conjecture}

\newcommand{\pf}{\noindent {\bf Proof.}\ }
\newcommand{\opf}{\noindent {\bf Outline of Proof.}\ }
\newcommand{\qed}{\hfill \hbox{\rule[-2pt]{3pt}{6pt}} \bigskip}

\newcommand{\mysection}[1]{\addtocounter{section}{1}%
\setcounter{thm}{0}\setcounter{equation}{0}%
\section*{\S{\Large \arabic{section}} #1}}%
\numberwithin{equation}{thm}

\newcounter{condsnum}
\numberwithin{condsnum}{thm}
\renewcommand{\thecondsnum}{(\thethm.\arabic{condsnum})}

\newcommand{\condsitem}{\item[{\em \thecondsnum}] \addtocounter{condsnum}{1}}

\newenvironment{conds}[1]{%
  \begin{list}{}
   {
   \settowidth{\labelwidth}{(1.1.11)}
   \setlength{\leftmargin}{1.1\labelwidth}}
  \setcounter{condsnum}{#1}
}%
{\end{list}}

\begin{document}
\setcounter{section}{-1} % The section starts from zero.

\title{Problems on geometric structures \\%
of projective embeddings}

\author{Takeshi Usa 
\thanks{
Dept. of Math., Himeji Institute of Technology, 
% \endgraf
2167 Syosya, Himeji, 671-2201 Japan.
\endgraf
{\it E-mail address  \/ }: usa{\char'100}sci.himeji-tech.ac.jp 
\endgraf
Typeset by \AmS-\LaTeX \ with P.Burchard's diagram package}}

\date{}

\maketitle

%============= Abstract =====================

\begin{abstract}
\noindent
This is a survey of our research 
on geometric structures of projective embeddings 
and includes some topics of our talks 
in several symposia during 1990-99. 
We clarify our main 
problem, which is to construct a kind of geometric 
composition series of projective embeddings. 
The concept of "geometric composition series" 
is an analogy in 
Algebraic Geometry with 
Jordan-H\"{o}lder series in Group Theory.
We present two of the candidates for 
the construction problem.
To approach this problem, 
we show several results and new tools for 
handling higher obstructions appearing 
in infinitesimal liftings. 
As a byproduct of the tools, 
we obtain a simplified proof for a criterion 
on arithmetic normality described 
in terms of Differential Geometry. \\

\end{abstract}

%%%%%%%%%%% Introduction %%%%%%%%%%%%%%%%%%%%%%

\mysection{Introduction.}
\noindent
In this article, we present several problems 
arising from our investigation 
during 1990-99 on geometric structures of 
projective embeddings (cf. \cite{PSE}, \cite{JPEXP} 
for partial reports).  When we use the technical term 
"geometric structure" of a projective embedding, 
it is our concern to see what kinds of intermediate 
ambient varieties appear for the projective 
subvariety defined by the given embedding.

To clarify this point more precisely, 
let us consider a connected complex projective manifold 
$X$ of dimension $n > 0$ and an embedding 
$j: X \hookrightarrow P={\Bbb P}^{N}({\Bbb C})$. 
Then, by an elementary fact on polynomial rings, 
we see that for any integer $q$ with $n < q < N$, 
there exists a projective subvariety $W$ of dimension $q$ 
satisfying $j(X) \subset W \subset P$. In this case, we 
say that the variety $W$ is an 
{\em intermediate ambient variety} 
of the subvariety $j(X)$.  

On the other hand, 
if we suppose an additional condition on $W$, 
e.g. a variety $W$ to be smooth 
along $j(X)$, namely $j(X) \subseteq Reg(W)$, 
then we can not assure the existence of a variety $W$ 
satisfying the condition.  For example, taking 
a Horrocks-Mumford abelian surface $A$ in 
$P={\Bbb P}^{4}({\Bbb C})$ as the subvariety $j(X)$, 
then there is no hypersurface $W$ with $A \subseteq Reg(W)$, 
which is certified by the calculation of Pic(W).

Thus we have special interest on the existence problem of 
intermediate ambient varieties with some additional 
conditions which can characterize the embedding. 
Then, we face an important problem, namely what 
conditions should be posed as the additional conditions?
One of the candidates for the condition 
is presented in Definition \ref{GSH}.
By using the concept "geometric shell", 
we can state our very optimistic 
Working Hypothesis \ref{SWHP}, 
which claims the existence 
of a projective embedding with 
a good decomposition by geometric shells. 
We should make a remark that 
this working hypothesis arose from the strong influence of 
the works done by Fujita, Mori, Mukai, and Sommese 
(e.g. \cite{CTPV}, \cite{OGCI}, \cite{CK3F3}, \cite{MCNAD}).
As an approach to this working hypothesis, we summarize in 
$\S 2$ the results on Lefschetz operators and on meta-Lefschetz 
operators. We also present Conjecture 
\ref{MCONJ} and 
clarify the relation with the former working hypothesis.
As a preparation for attacking the conjectures, 
we newly introduce several key concepts 
for the infinitesimal method in $\S 3$.
They often help us to remove the difficulty of higher 
obstructions for making 
the correspondence between subsheaves of the normal bundle 
and intermediate ambient varieties. 
In $\S 4$, we discuss 
arithmetic normality 
from two points of view. 
The first view point concerns our 
framework and strategy for studying geometric 
structures of projective embeddings. 
From the second point of view, namely that 
of Differential Geometry, we explain 
a criterion for arithmetic normality 
in terms of the second fundamental form. 
Here we describe an outline of 
another proof for the criterion 
which is simplified by the tools in $\S 3$.
This will show the power of our new tools.

In this article, we consider only 
the objects defined over the complex 
number field ${\Bbb C}$. In case of handling 
graded objects, we consider only homomorphisms  
of preserving their grading otherwise mentioned. 
For example, "minimal free resolution" always means 
"graded minimal free resolution". 
Sometimes we state our results by using a pair 
of varieties and a pair of their embeddings 
instead of using the term "subvarieties". That is only to 
emphasize the fact that we can choose the embeddings 
suitably with fixing the pair of varieties in the real 
situation.

The author deeply thanks to Prof.O.A.Laudal and Prof.S.J.Kwak 
for their warmful encouragement, 
to Prof.M.Hashimoto for showing me a nice fact, 
and to Prof.S.Tsuboi and 
Prof.C.Miyazaki for their heavy efforts of organizing symposia 
where one could have precious chances to meet the former 
two people.

%%%%%%%%%  Working Hypothesis  %%%%%%%%%%%%%%%%%%%

\mysection{Working Hypothesis.}
\noindent
In this section,  we present a 
key concept for considering 
geometric structures of embeddings 
and show several problems, in particular 
our optimistic working hypothesis. 
We hope that this may also bring us an insight 
for studying the syzygies of projective 
subvarieties.

\bigskip

Let us confirm our notation used 
in the sequel.

%===== Notation =======

\begin{noten} \label{NTA}
Let us take a complex projective scheme $X$ of dimension $n$ 
and one of its embeddings
$j: X \hookrightarrow P={\Bbb P}^{N}({\Bbb C})$.
The sheaf of ideals defining $j(X)$ in $P$ and
the conormal sheaf are denoted 
by $I_{X}$ and $N^{\vee}_{X/P}=I_{X}/I_{X}^2$,
respectively.  Taking a ${\Bbb C}$-basis 
$\{ Z_0, \ldots ,Z_N \}$ of $H^0(P,O_P(1))$. 
Then we put 

\[
\begin{array}{rcl}
S &:=& \bigoplus\limits _{m \geq 0} H^0(P,O_P(m)) %
\cong {\Bbb C}[Z_0, \ldots ,Z_N]\\
S_{+} &:=& =\bigoplus\limits _{m > 0} H^0(P,O_P(m)) %
\cong (Z_0, \ldots ,Z_N){\Bbb C}[Z_0, \ldots ,Z_N] \\
\widetilde{R_X} &:=& \bigoplus\limits_{m \geq 0} H^0(X,O_X(m)) \\
{\Bbb I}_X &:=& \bigoplus\limits_{m \geq 0} H^0(P, I_X(m)) \\
R_X &:=& Im[S \rightarrow \widetilde{R_X}] \cong S/{\Bbb I}_X  \\
gsyz_X^q(m) &:=& Tor^S_q(R_X,S/S_{+})_{(m)},
\end{array}
\]

\noindent
where the subscript $(m)$ of $Tor$ above means taking 
its degree $m$ part as a graded $S$-module. 
Obviously, the space $gsyz_X^q(m)$ represents 
minimal generators in degree $m$ of the 
q-th syzygy of $R_X$ as an $S$-module.

\end{noten}

\bigskip

As a preparation, we recall the following key 
concepts introduced in \cite{EML}.

%====== Defn. of G-shell ==========

\begin{defn}[PG-shell and G-shell] \label{GSH}
Let $V$ and $W$ be closed subschemes of 
$P={\Bbb P}^{N}({\Bbb C})$ which satisfy $V \subseteq W$ 
(In this case, the subscheme $W$ is called simply 
an {\em intermediate ambient scheme} 
of $V$). If the natural map:

\[
\mu_q : Tor^S_q(R_W,S/S_{+}) \rightarrow Tor^S_q(R_V,S/S_{+}) 
\]

\noindent
is injective for every $q \geq 1$, we say that $W$ is 
a {\em pregeometric shell (abv. PG-shell)} of $V$. Moreover, 
if $W$ is a closed subvariety and the 
regular locus $Reg(W)$ of $W$ contains $V$, we say that 
$W$ is a {\em geometric shell (abv. G-shell)} of $V$.
For the subscheme $V$, $P$ and $V$ itself are called 
{\em trivial PG-shell (or trivial G-shell)}.
\end{defn}

\bigskip

Now let us see several elementary 
facts relating with "PG-shell".

%======== Prop of PG-shells ============

\begin{prop}
Let $V$ and $W$ be closed subschemes of 
$P={\Bbb P}^{N}({\Bbb C})$ which satisfy $V \subseteq W$. 

\begin{conds}{1}

\condsitem If $W$ is a hypersurface, then $W$ is 
a PG-shell of $V$ if and only if the equation 
$H_W$ of $W$ is a member of minimal generators 
of the homogeneous ideal ${\Bbb I}_V$ of $V$.  

\condsitem Assume that the subscheme $V$ is 
a complete intersection. Then the scheme 
$W$ is a PG-shell of $V$ if and only if 
the subscheme $W$ is defined by a part of 
minimal generators of ${\Bbb I}_V$.

\condsitem Take a closed scheme $Y$ such that 
$V \subseteq Y \subseteq W$. Assume that $W$ is 
a PG-shell of $V$. Then $W$ is also a PG-shell of 
$Y$. In particular, the subscheme $W$ is also a PG-shell 
of the $m$-th infinitesimal 
neighborhood $Y=(V/W)_{(m)}$ 
of $V$ in $W$, where $(V/W)_{(m)}=(|V|,O_W/I_V^{m+1})$.

\condsitem Fix the subscheme $V$ of $codim(V,P) \geq 2$. 
Then all 
non-trivial PG-shells of $V$ form non empty algebraic 
family of finite components (N.B. The family of all 
non-trivial G-shells of $V$ may be empty even if 
$V$ itself is a smooth variety). 

\condsitem If $W$ is a PG-shell of $V$, then 
we have an inequality on their arithmetic depth: 
$arith.depth(V) \leq arith.depth(W)$. In particular, 
if the natural restriction map 
$H^0(P,O_P(m)) \rightarrow H^0(V, O_V(m))$ is 
surjective for all integers $m$
(i.e. $R_V = {\widetilde{R_V}}$), then 
the natural restriction map 
 $H^0(P,O_P(m)) \rightarrow H^0(W, O_W(m))$ is 
also surjective for all integers $m$
(i.e.$R_W = {\widetilde{R_W}}$). 

\condsitem If the subscheme $W$ is a PG-shell of 
the subscheme $V$ with $arith.depth(V) \geq 2$, 
then we have an inequality 
on their Castelnuovo-Mumford regularity: 
$CMreg(V) \geq CMreg(W)$.

\condsitem Assume that there exist $r$ hypersurfaces 
$H_1, \ldots H_r$ in $P$ which form $O_W$-regular 
sequence and satisfy $V=W \cap H_1 \cap \ldots \cap H_r$.  
If the restriction map 
$H^0(P,O_P(m)) \rightarrow H^0(V, O_V(m))$ is 
surjective for all integers $m$, then $W$ is 
a PG-shell. 

\condsitem  Assume that the subscheme $V$ is 
non-degenerate, namely no hyperplane contains 
$V$. If $W$ has a 2-linear resolution, i.e. 
the homogeneous coordinate ring $R_W$ of $W$ 
has a minimal $S$-free resolution of the form : 
$ 0 \leftarrow R_W \leftarrow S \leftarrow %
F_1(-2) \leftarrow F_2(-3) \leftarrow %
\cdots \leftarrow F_p(-p-1) \leftarrow \cdots $, 
where $F_u(v)$ denotes $\oplus S(v)$ : a direct sum of 
several copies of $S$ with 
degree $v$ shift, then $W$ is a PG-shell of $V$. 

\end{conds}
\end{prop}

\pf 
Directly from Definition \ref{GSH}, we see 
(1.3.1),(1.3.2) and (1.3.3). 
To get (1.3.4), using the injectivity of $Tor^S_1$, 
we notice that any PG-shell of $V$ is defined by 
a part of minimal generators of the ideal ${\Bbb I}_V$. 
Let us take a parametrizing space $T$ of all the intermediate 
ambient schemes of $V$ which are defined 
scheme-theoretically by parts 
of the minimal generators of the ideal ${\Bbb I}_V$. 
Obviously the parametrizing space $T$ is identified 
with a set-theoretic direct sum of Zariski open sets in 
several products of Grassmannian varieties. 
Through flattening stratifications including vertices 
of the affine cones of the members in the family, 
we get an algebraic family 
of intermediate ambient schemes of $V$ with constant
Betti numbers, whose parametrizing space is named $T$ again. 
We may assume that every component in $T$ 
includes a point for a PG-shell of $V$.
This family includes all the PG-shells of $V$. 
Looking at an induced chain homomorphism of (relatively) 
minimal free resolutions of ${\Bbb I}_V$ and of the family, 
we have only to extract the open sets of which corresponds 
to the "maximal" rank locus of every map in the chain homomorphism.  
For (1.3.6), use the Eisenbud-Goto criterion on 
Castelnuovo-Mumford regularity in \cite{EG}. Similarly (1.3.5) 
is obtained by applying the formula on depth and homological 
dimension. (1.3.7) is shown in \cite{EML}. On the claim (1.3.8), 
see Lemma1 in \cite{FRCA} from our point of view.
\qed

\bigskip

%============  Hashimoto =======================

The next fact is kindly told me by Prof.M.Hashimoto 
with answering some questions relating the concept 
of G-shells. It may help us to construct 
G-shells in the real situation.

\begin{prop}
Let $Y$ and $Z$ be arithmetically Cohen-Macaulay 
projective subschemes of $P={\Bbb P}^N({\Bbb C})$. 
Assume that $X=Y \cap Z$ is of codimension $a+b$, 
where $a=codim(Y,P)$ and $b=codim(Z,P)$. 
Then both $Y$ and $Z$ are PG-shells of $X$.
\end{prop}

\pf
For $R_Y$ and $R_Z$, take their minimal $S$-free 
resolutions: ${\Bbb F}_{\bullet} \rightarrow R_Y$, 
${\Bbb G}_{\bullet} \rightarrow R_Z$, whose length are 
$a$ and $b$, respectively. If we can show that 
the complex 
${\Bbb F}_{\bullet} \otimes {\Bbb G}_{\bullet}$ 
is acyclic, then triviality of the complex 
$({\Bbb F}_{\bullet} \otimes {\Bbb G}_{\bullet}) %
\otimes (S/S_{+}) \cong 
({\Bbb F}_{\bullet} \otimes S/S_{+})\otimes %
({\Bbb G}_{\bullet}\otimes S/S_{+})$ 
means that 
the complex 
${\Bbb F}_{\bullet} \otimes {\Bbb G}_{\bullet}$ 
is a minimal $S$-free resolution of $R_{X} \cong R_Y \otimes R_Z$. 
Since the complexes ${\Bbb F}_{\bullet}$ and ${\Bbb G}_{\bullet}$ 
are naturally considered as subcomplexes and as direct summands 
of the complex 
${\Bbb F}_{\bullet} \otimes {\Bbb G}_{\bullet}$, we see that 
the schemes $Y$ and $Z$ are PG-shells of $X$.
The complex 
${\Bbb F}_{\bullet} \otimes {\Bbb G}_{\bullet}$ 
has the length $a+b$, 
which coincides with 
$ht({\Bbb I}_Y+{\Bbb I}_Z)=depth({\Bbb I}_Y+{\Bbb I}_Z,S)$.
To see the acyclicity of 
${\Bbb F}_{\bullet} \otimes {\Bbb G}_{\bullet}$, 
we apply 
Buchsbaum-Eisenbud criterion for acyclicity on free complexes 
(cf. \cite{WMCE}). 
Thus we have only to show that for any prime ideal 
${\frak p} \in Spec(S)$ with $depth({\frak p}) < a+b$, 
$({\Bbb F}_{\bullet} \otimes {\Bbb G}_{\bullet})_{\frak p}$ 
is acyclic. If $ht({\frak p})=depth({\frak p}) < a+b$, then 
$\frak p \not\supset {\Bbb I}_Y+{\Bbb I}_Z$, namely 
$\frak p \not\supset {\Bbb I}_Y$ or 
$\frak p \not\supset {\Bbb I}_Z$. For example, if 
$\frak p \not\supset {\Bbb I}_Y$ , then 
$({\Bbb F}_{\bullet})_{\frak p} \rightarrow 0$ is 
split exact and therefore 
$({\Bbb F}_{\bullet} \otimes {\Bbb G}_{\bullet})_{\frak p}$ 
is acyclic.
\qed

The following example shows that all the exceptional 
cases in the classical Petri's Analysis 
can be considered as the cases of G-shells appearing.

%======= Example of  Petri's Analysis =========

\begin{exam}[Quadric hulls in Petri's Analysis]
Let $X=C$ be a non-hyperelliptic smooth projective 
curve of genus $g \geq 4$, and 
$j=\Phi_{|K_C|}:C \hookrightarrow P={\Bbb P}^{g-1}({\Bbb C})$
its canonical embedding.
Taking the quadric hull $W$ of $j(C)$, 
namely the closed subscheme defined by all equations 
of $j(C)$ with degree $2$. Then, 
the quadric hull $W$ coincides with $j(C)$ itself or 
is a non-trivial G-shell of $j(C)$. 
\end{exam}

\pf
If $g=4$, then the surface $W$ may have a singular point 
but the embedded curve $j(C)$ is a non-singular 
complete intersection of type $(2,3)$. Thus we may assume 
that $g \geq 5$.  By classical 
Petri's Analysis (cf. \cite{DAFV}, \cite{PALSQ}), we see 
the exceptional cases explicitly, namely $W$ is a Veronese 
surface in ${\Bbb P}^5({\Bbb C})$ or a rational normal scroll. 
In both cases, $W$ is a surface of minimal degree. 
Then apply (5.2)Lemma in \cite{SCC} (see also \cite{ERS}), 
we obtain that $W$ has 2-linear resolutions, which implies 
that $W$ is a G-shell of $j(C)$.
\qed

\bigskip

%======== Problems 0n G-shell =====================

\begin{prob}
To make a foundation for studying PG-shells or G-shells, 
let us list several problems conjured up naturally 
in our mind. 

\begin{conds}{1}

\condsitem For a non-hyperelliptic 
curve $C$ of genus $g=g(C) \geq 3$
and its canonical embedding 
$j=\Phi_{|K_C|}:C \hookrightarrow P={\Bbb P}^{g-1}({\Bbb C})$, 
classify all the  PG-shells of $j(C)$ 
{\em (cf. Green Conjecture \cite{KGPV}, \cite{FRCA})}.

\condsitem Describe the condition of "PG-shell" in 
terms of "generic initial ideals".

\condsitem  Assume that a projective subscheme $W$ is a PG-shell 
of a projective subvariety $V \subset P={\Bbb P}^N({\Bbb C})$. 
Then the subscheme $W$ is always reduced and irreducible ? (N.B. When 
the subscheme $W$ is a hypersurface, this is true.)

\condsitem Take smooth projective subvarieties 
$V$ and $W$ of positive dimension.  
Assume that the subvariety $V$ is arithmetically normal. 
If $W$ is a PG-shell of $V$, then 
does the inequality on $\Delta$-genus {\em (cf. \cite{CTPV})}:
$\Delta (V,O_P |_V (1)) \geq \Delta (W,O_P |_W (1))$ 
hold in general? (e.g. As a typical case, if the polarized 
manifold $(V,O_P |_V (1))$ is a member 
of a ladder of the polarized manifold $(W,O_P |_W (1))$, 
then $W$ is a G-shell of $V$ and 
this inequality is obviously true. 
On the other hand, if we assume that 
$V$ is non-degenerate arithmetically Buchsbaum, 
and $W$ is a hypersurface, then 
the result \cite{CBLCMS} on special cases of 
Eisenbud-Goto conjecture shows that this claim 
is also true.)

\condsitem Take a smooth projective 
subvariety $V$, a vector bundle $E$ on $V$, 
a section $\sigma \in \Gamma (V,E)$ 
which is transverse to the zero section, 
and its zero locus $X=Z(\sigma)$. Assume that 
$V$ is a G-shell of $X$. Then is the 
bundle $E$ always nef ?

\condsitem Take a smooth projective 
subvariety $V \subset P={\Bbb P}^N({\Bbb C})$ 
of dimension $n \geq 5 $. 
Assume that $V$ is not a hypersurface
and has no non-trivial G-shell. 
Then $codim(V) \geq n/2$? 
{\em (Implied by Hartshorne's C.I.conjecture. cf. \cite{VSC} 
\cite{TSAV})}
Moreover, for any positive integer $M$ and a linear 
embedding $P={\Bbb P}^N({\Bbb C}) \subset Q={\Bbb P}^{N+M}({\Bbb C})$, 
if the subvariety $V$ has no non-trivial G-shell except 
(multiple) projective 
cones, then does the Kodaira dimension of $V$ satisfy 
the inequality $\kappa (V) < 0$?

\end{conds}
\end{prob}

\bigskip

Now we present our working hypothesis in 
the most optimistic version, which 
suggests the direction of our 
research aiming.

%====== Working Hypothesis =======

\begin{whypo} \label{SWHP}
Let $X$ be a connected complex projective manifold of 
dimension $n > 0$. Then there exists an embedding: 
$j: X \hookrightarrow P={\Bbb P}^{N}({\Bbb C})$, which 
satisfies the following conditions.

\begin{conds}{1}
\condsitem There is a set of G-shells 
$ \{ W_p \} _{p=0}^k$ of $j(X)$ 
which satisfy:
$j(X)=W_0 \subset W_1 \subset \ldots \subset W_k \subseteq P $ 
and moreover $W_{p-1} \subset Reg(W_{p})$ for $p=1,\ldots ,k$.

\condsitem For each $p=1,\ldots ,k$, 
there is a "nef" vector bundles $E_p$ on $W_p$ and a section 
$\sigma_p \in \Gamma (W_p,E_p)$ such that the zero locus 
$Z(\sigma_p)$ coinsides with $W_{p-1}$ and 
$rank(E_p)=dim(W_p)-dim(W_{p-1})$.

\condsitem The subvariety $W_k$ 
has a birational morphism from 
a projective bundle over 
a homogeneous space (in the sense of including abelian varieties). 
\end{conds}

\noindent
The set $\Xi=\{ (W_p,E_p,\sigma_p) \} _{p=1}^k$ of $j(X)$ 
and the integer $k$ 
are called {\em a geometric composition series} of the embedding j 
or of the subvariety $j(X)$ and {\em the length } of
the geometric composition series $\Xi$, respectively.
For a given projective 
manifold $X$, if the embedding $j_0$ has a geometric composition 
series $\Xi_0$ whose length $k_0$ attains the minimum among 
the embeddings of $X$ with geometric composition series, 
then we say that the geometric composition series $\Xi_0$ 
is {\em a absolutely minimal geometric composition series} of $X$.   
\end{whypo}

\bigskip

%===== Remark on WHypo =========

\begin{remk} 

To avoid confusion or to clarify 
what is in the author's 
mind, one should describe several 
points.

\begin{conds}{1}

\condsitem For a vector bundle $E$ on a 
projective variety $V$,  
we say that the bundle $E$ is {\em nef} 
if the tautological 
line bundle $L_E=O_{P(E)/V}(1)$ is nef 
on the projective 
bundle $P(E)={\Bbb P}(E)$ over $V$ 
associated to 
the bundle $E$, namely for any 
curve $C$ in $P(E)$, the intersection 
number satisfies the inequality:
$(L_E.C) \geq 0$. 

\condsitem Frankly speaking, 
the author confess that 
we might have to weaken 
our working hypothesis 
to some extent 
in the real situation.
For example, we might have to 
replace the conditions: 
{\em (a)} "PG-shells" instead of "G-shells" ; 
{\em (b)} "reflexive sheaves" in stead of 
"vector bundles" ;
{\em (c)}"rather mild singular locus of $W_p$" 
instead of "$Reg(W_{p})$" ; 
{\em (d)} "$\kappa (W_k) \leq 0$" instead of     
"a homogeneous space." 

\end{conds}
\end{remk}

\bigskip

%========Comp. Ser. & Cpl.Int.=================

\begin{prop}
Let $X$ be a connected complex projective manifold of 
dimension $n \geq 2$ and 
$j: X \hookrightarrow P={\Bbb P}^{N}({\Bbb C})$
an embedding. Then the following four conditions 
are equivalent.

\begin{conds}{1}

\condsitem The subvariety $j(X)$ is a complete intersection.

\condsitem There is a set of intermediate 
ambient varieties  
$ \{ W_p \} _{p=0}^{N-n}$ of $j(X)$ 
which satisfies the conditions:
{\em (a)}$dim(W_p)=n+p$ ;
{\em (b)} $j(X)=W_0 \subset W_1 \subset \ldots %
\subset W_{N-n}= P $ ;
{\em (c)} $W_{p-1} \subset Reg(W_{p})$ 
for $p=1, \ldots , N-n$.

\condsitem The embedding $j$ has a 
geometric composition series 
$\Xi= \{ (W_p,E_p,\sigma_p) \}_{p=1}^{N-n}$ of 
length $N-n$ with $rank(E_p)=1$. 

\condsitem The embedding $j$ has a 
%%% minimal 
geometric composition series $\Xi_0= \{ (W_1,E_1,\sigma_1) \}$ of 
length $1$ which satisfies $W_1=P$ and 
$E_1=\oplus_{s=1}^{N-n} O_P(m_s)$. 

\end{conds}
\end{prop}

\pf
The essential part is to show the equivalence 
between (1.9.1) and (1.9.2). Assume that 
(1.9.2). Starting from $W_{N-n}$ and using that 
each $W_p$ is a Cartier divisor of $W_{p+1}$, 
we show inductively 
that each $W_p$ is a complete intersection 
and $Pic(W_p) \cong {\Bbb Z}O_{W_p}(1)$ 
for $p \geq 1$ by virtue of 
Corollary 3.2 in \cite{ASAV}, which is still 
valid in the singular cases. Thus we have (1.9.1).
Contrary, now we assume (1.9.1). 
A little care is needed to 
apply Bertini's theorem and 
to see that $W_{p-1} \subset Reg(W_{p})$, 
which is rather a strong condition than 
$X \subset Reg(W_{p})$.
Take hypersurfaces $D_1, \ldots , D_r$ of degree 
$d_1 , \ldots , d_r$,
respectively such that $r=N-n$, 
$j(X)= D_1 \cap \ldots \cap D_r$ 
and $d_1 \leq \ldots  \leq d_r$. 
Then consider the linear system 
$\Lambda_r = H^0(P,I_X(d_r))$ on $P=W_{r}$. 
Since $I_X(d_r)$ is generated by global sections,
the base locus $Bs(\Lambda_r)$ coincides with $X$. 
Also by $D_r \in \Lambda_r$ satisfying $X \subset Reg(D_r)$, 
we find that general members are smooth. Then we 
put $W_{r-1}$ to be a smooth member of $\Lambda_r$. 
Obviously 
$j(X)=D_1 \cap \ldots \cap D_{r-1} \cap W_{r-1}$. 
As an induction 
hypothesis, we may assume that we have 
smooth complete intersection subvarieties: 
$W_{k}, W_{k+1},\ldots, W_r=P$ such that 
$dim(W_p)=n+p$, 
$j(X)=D_1 \cap \ldots \cap D_p \cap W_{p}$
for $p=k, \ldots ,r$.  We may assume $k \geq 2$.
Then we consider a sublinear system 
$\Lambda_k :=H^0(W_{k},I_{X/W_{k}}(d_k)) \subset H^0(W_{k},O_{W_{k}}(d_k))$ 
on the subvariety $W_k$. Since $I_{X/W_{k}}(d_k)$ is generated by 
the sections $D_1, \ldots, D_k$, 
namely 
$(D_1, \ldots, D_k):%
\oplus_{q=1}^k O_{W_k}(d_k-d_q) %
\rightarrow %
I_{X/W_{k}}(d_k)$ is 
surjective, we have $Bs(\Lambda_k)=X$. 
By the same argument as above, we obtain a smooth member 
$W_{k-1} \in \Lambda_k$. Then, using the arithmetic normality 
of $W_k$, it is easy to see that $W_{k-1}$ 
is also a complete intersection and 
$j(X)=D_1 \cap \ldots \cap D_{k-1} \cap W_{k-1}$. 
\qed

%%%%%%%%%%  Conjectures  %%%%%%%%%%%%%%%%%%%%%

\mysection{Conjectures.}
\noindent
In this section, we give some conjectures 
relating to Lefschetz operators. We expect 
that these conjectures give an approach 
to get our previous working hypothesis.  

First, let us recall the definition of 
Lefschetz operators (cf. \cite{LEO}).

%======== Defn. usual Lef. op. ==========

\begin{defn}[Lefschetz operator]
Let $X$ be a complex projective scheme of dimension $n \geq 0$, 
$j: X \hookrightarrow P={\Bbb P}^{N}({\Bbb C})$ an embedding, 
$E$ an ${\cal O}_X$-coherent sheaf, and 
$N^{\vee}_{X/P}$ the conormal sheaf of $j(X)$ in $P$, 
where $I_{X}$ denotes the sheaf of ideals defining $j(X)$ in $P$.
By natural restriction:
$j^{*}:H^1(P,\Omega_P^1) \rightarrow H^1(X,\Omega_X^1)$, 
we have a hyperplane class 
$h=j^*(c_1(O_P(1))) \in H^1(X,\Omega_X^1)$, which induces 
a cohomological operator (depending on the embedding $j$): 

\[
\begin{CD}
L_X:H^p(X,\Omega_X^q \otimes E) @>{h \cup}>> %
H^{p+1}(X,\Omega_X^{q+1}\otimes E)
\end{CD}
\]

\noindent
For a section $\sigma \in H^0(X,E)$, if the class 
$L_X^p(\sigma) \in H^p(X,\Omega_X^p \otimes E)$ is 
not zero and $L_X^{p+1}(\sigma)$ is zero, then we 
say that the section $\sigma$ has {\em the penetration order} 
$p$ and denote it by $pent(\sigma)=p$. For an equation 
$F \in H^0(P,I_X(m))$ of $j(X)$ with degree $m$, we define 
$pent(F)=pent([F])$ by putting $E=N^{\vee}_{X/P}(m)$, 
where $[F]$ denotes the section of 
$H^0(X,N^{\vee}_{X/P}(m))$ induced by the natural restriction 
$I_{X} \rightarrow I_{X}/I_{X}^2=N^{\vee}_{X/P}$.
\end{defn}

\bigskip

We introduce meta-Lefschetz operators, which 
are difficult to control but give finer 
information than Lefschetz operators.

%=================== Defn meta-Lef. op. ===============988-1115

\begin{defn}[meta-Lefschetz operator \cite{EML}, \cite{MLO}]
Let $X$ and $W$ be complex projective schemes of dimension 
$n \geq 0$ and of dimension $N$, respectively.  
Take an arbitrary line bundle $O_W(1)$ on $W$ 
and an embedding $j: X \hookrightarrow W$. 
Assume that $j(X) \subset Reg(W)$.  Put 
$N^{\vee}_{X/W}=I_{X}/I_{X}^2$ to be 
the conormal sheaf of $j(X)$ in $W$, 
where $I_{X}$ denotes the sheaf of 
ideals defining $j(X)$ in $W$. 
Then we take the de Rham complex 
$\Omega_W^{\bullet}$ of $W$ : \\

\[
\begin{CD}
0 @>>> O_W @>d>> \Omega_W^1 @>d>> \Omega_W^2 @>d>> {\ldots} @>d>> %
\Omega_W^N @>>> 0 \\
\end{CD}
\]

\noindent
and the ideal order filtration {\em (cf. \cite{FHF})} 
$F^p_{\nu}\Omega_W^{\bullet}$:

\[
\begin{CD}
0  @>>> I_{X}^{\nu + p} @>d>> I_{X}^{\nu + p-1} \Omega_W^1 %
@>d>> {\ldots} @>d>> \\ 
    I_{X}^{\nu + 1} \Omega_W^{p-1}  @>d>> \Omega_W^{p} @>d>> %
               {\ldots} @>d>> \Omega_W^N @>>> 0. \\
\end{CD}
\]

\noindent
Now we fix $\nu$ and see 
$Gr^{p}_{F^{\nu}}(\Omega_W^{\bullet})=F^p_{\nu}/F^{p+1}_{\nu}$:

\[
\begin{CD}
0  @>>> I_{X}^{\nu + p}/I_{X}^{\nu + p+1} @>\overline{d_I}>> %
I_{X}^{\nu + p-1}/I_{X}^{\nu + p} \otimes \Omega_W^1 %
@>\overline{d_I}>> {\ldots} \\ 
@>\overline{d_I}>> I_{X}^{\nu + 1}/I_{X}^{\nu + 2} %
\otimes \Omega_W^{p-1}  %
@>\overline{d_I}>> \Omega_W^{p}|_{X_{(\nu)}} @>>> 0,  \\
\end{CD}
\]

\noindent
where $X_{(\nu)} = (|X|,O_{W}/I_X^{\nu+1})$. Contrary 
to the fact that
the exterior derivative $d$ is not $O_W$-linear, the map 
$\overline{d_I}$ is $O_{W}$-linear and  
compatible with tensoring 
by $O_W(m)$. Thus we have:

\[
\begin{CD}
I_{X}^{\nu + 1}/I_{X}^{\nu + 2}(m) \otimes \Omega_W^{p-1}  %
@>\overline{d_I}>> \Omega_W^{p}|_{X_{(\nu)}}(m)
\end{CD}
\]

\noindent
and 

\[
\begin{CD}
H^t(X, I_{X}^{\nu + 1}/I_{X}^{\nu + 2}(m) \otimes \Omega_W^{p-1}) %
@>\overline{d_I}>>H^t(X_{(\nu)}, \Omega_W^{p}|_{X_{(\nu)}}(m))
\end{CD}
\]

\noindent
Next we consider a natural exact sequence ($\overline{LFT}$):

\[
\begin{CD}
0 @>>> I_{X}^{\nu + 1}/I_{X}^{\nu + 2} \otimes \Omega_W^{p}(m)%
 @>>>\Omega_W^{p}(m)|_{X_{(\nu+1)}} @>>>%
\Omega_W^{p}(m)|_{X_{(\nu)}} @>>> 0,
\end{CD}
\]

\noindent
which induces an obstruction map:

\[
\overline{\delta}_{LFT}^{(\nu)}:H^s(X_{(\nu)}, %
\Omega_W^{p}(m)|_{X_{(\nu)}}) %
\rightarrow H^{s+1}(X, %
I_{X}^{\nu + 1}/I_{X}^{\nu + 2} \otimes \Omega_W^{p}(m)).
\]

\noindent
Then we can define a map:

\[
\widehat{L}_{X/W}^{(\nu)}=\overline{\delta}_{LFT}^{(\nu)} %
\circ \overline{d_I}:H^{a}(X, %
I_{X}^{\nu + 1}/I_{X}^{\nu + 2}(m) \otimes \Omega_W^{b})
\rightarrow %
H^{a+1}(X, %
I_{X}^{\nu + 1}/I_{X}^{\nu + 2}(m) \otimes \Omega_W^{b+1}),
\]

\noindent
which is called the {\em $\nu$-th meta-Lefschetz operator} with 
respect to the embedding $j:X \hookrightarrow W$.
In case that $W=P={\Bbb P}^{N}({\Bbb C})$ and 
$O_P(1)$  is the tautological ample line bundle, 
the symbol $\widehat{L}_{X/W}^{(\nu)}$ is 
abbreviated to $\widehat{L}_{X}^{(\nu)}$.  
Moreover, if $\nu = 0$, we denote it by $\widehat{L}_X$ 
instead of $\widehat{L}_X^{(0)}$ and call it simply
{\em meta-Lefschetz operator} if there is no danger 
of confusion. 
For the meta-Lefschetz operator $\widehat{L}_X$, 
we set 

\[
\widehat{gsyz}_X^q(m):=Im[\widehat{L}_X :%
H^{0}(X, %
N^{\vee}_{X/W}(m) \otimes \Omega_W^{q-1})
\rightarrow %
H^{1}(X, %
N^{\vee}_{X/W}(m) \otimes \Omega_W^{q})]
\]

\end{defn}

\bigskip

Fundamental properties on meta-Lefschetz operator 
are given as follows.

%============ Thm on meta-Lef.op. ============= 

\begin{thm}[\cite{MLO}] \label{MLTH}
Let $X$ be a complex projective variety of dimension 
$n$, 
$j: X \hookrightarrow P={\Bbb P}^{N}({\Bbb C})$ an embedding, 
$N^{\vee}_{X/P}=I_{X}/I_{X}^2$ the conormal sheaf of $j(X)$ in $P$, 
where $I_{X}$ denotes the sheaf of ideals defining $j(X)$ in $P$.
Take the meta-Lefschetz operator $\widehat{L}_X$ with respect to 
the embedding. Then the following properties hold.

\begin{conds}{1}

\condsitem The meta-Lefschetz operator has naturality. 
In other words, for any closed subscheme Y, the diagram:

\[
\begin{CD}
H^{p}(X, %
N^{\vee}_{X/P}(m) \otimes \Omega_P^{q})
@>{\widehat{L}_X}>> %
H^{p+1}(X, %
N^{\vee}_{X/P}(m) \otimes \Omega_P^{q+1}) \\
@V{natural}VV   @VV{natural}V \\
H^{p}(Y, %
N^{\vee}_{Y/P}(m) \otimes \Omega_P^{q})
@>>{\widehat{L}_Y}> %
H^{p+1}(Y, %
N^{\vee}_{Y/P}(m) \otimes \Omega_P^{q+1}) 
\end{CD}
\]

\noindent
is commutative.

\condsitem The diagram:

\[
\begin{CD}
H^{p}(X, %
N^{\vee}_{X/P}(m) \otimes \Omega_P^{q})
@>{\widehat{L}_X}>> %
H^{p+1}(X, %
N^{\vee}_{X/P}(m) \otimes \Omega_P^{q+1}) \\
@V{natural}VV   @VV{natural}V \\
H^{p}(X, %
N^{\vee}_{X/P}(m) \otimes \Omega_X^{q})
@>>{-m \cdot L_X}> %
H^{p+1}(X, %
N^{\vee}_{X/P}(m) \otimes \Omega_X^{q+1}) 
\end{CD}
\]

\noindent
is commutative, where $L_X$ denotes the Lefschetz operator.

\condsitem Assume that $j(X)$ has arithmetic depth 
$\geq 2$, which includes the case that
$X$ is a normal projective
variety of dimension $n > 0$ and the embedding is 
arithmetically normal,
namely $H^0(P,O_P(m)) \rightarrow H^0(X, O_X(m))$ is surjective 
for all integers $m$. Then there is a natural one to one 
correspondence 
$\gamma^q(m): gsyz_X^q(m) \rightarrow \widehat{gsyz}_X^q(m)$ 
as vector spaces. Here the space 
$ gsyz_X^q(m)$ represents minimal 
generators in degree $m$ of the q-th syzygy of $R_X$.  

\condsitem  For an integer $k$ satisfying $n-1 \geq k \geq 1$, 
assume that 
the projective subvariety $j(X)$ 
has arithmetic depth 
$k+2$, or equivalently $H^s(X, O_X(u))=0$ for 
$u \in {\Bbb Z}$, $k \geq s \geq 1$, 
and $H^0(P,O_P(m)) \rightarrow H^0(X, O_X(m))$ is surjective 
for all integers $m$. 
Then the $k$-uple of the meta-Lefschetz operator:

\[
(\widehat{L}_X)^k :%
H^{1}(X, %
N^{\vee}_{X/P}(m) \otimes \Omega_P^{q})
\rightarrow %
H^{k+1}(X, %
N^{\vee}_{X/P}(m) \otimes \Omega_P^{q+k})
\]

\noindent
is injective on the subspace $\widehat{gsyz}_X^q(m)$ for
all integers $m$. 
Moreover, the map 
$\overline{d_I}:H^{k+1}(X,N^{\vee}_{X/P}(m) \otimes \Omega_P^{q+k})%
\rightarrow H^{k+1}(X,\Omega_P^{q+k+1} \otimes O_X(m))$ 
is injective on the subspace 
$(\widehat{L}_X)^k(\widehat{gsyz}_X^q(m))$.
\end{conds}
\end{thm}

\bigskip

Returning to Lefschetz operators and 
make a preparation for defining 
Lefschetz chains and dual Lefschetz chain 
which play key roles in our conjectures.

\bigskip

%======== Making a Lefschetz filtration ============

Now we take the canonical map 
$\gamma^1(m): gsyz_X^1(m) \rightarrow \widehat{gsyz}_X^1(m)$ 
in the Theorem\ref{MLTH} above for $q=1$ and consider the 
commutative diagram:

\[
\begin{CD}
\Sigma Z_t H^0(I_X(m-1)) %
@>>> H^0(I_X(m)) @>>> gsyz^1_X(m) @>>> 0 \\
@. @V{natural}VV  @VV{\gamma^1(m) \cong}V \\
@. H^0(N^{\vee}_X(m)) @>>{\widehat{L}_X}> \widehat{gsyz}^1_X(m) %
@>>> 0  \\
@. @V{L_X}VV  @VV{Inclusion}V @. \\
@. H^1(N^{\vee}_X(m)\otimes \Omega^1_X) @<<{natural}< %
H^1(N^{\vee}_X(m)\otimes \Omega^1_P|_X) @. \\
@. @V{L_X^{p-1}}VV @. \\
@. H^p(N^{\vee}_X(m)\otimes \Omega^p_X),  @. @.
\end{CD}
\]

\noindent
where the first and the second rows are exact. Then we put:

\[
J_p(m)=Ker[H^0(I_X(m)) \rightarrow H^0(N^{\vee}_X(m)) %
\stackrel{L_X^{p+1}}{\rightarrow}H^{p+1}(N^{\vee}_X(m)\otimes \Omega^{p+1}_X)].
\]

\noindent
From a chain of ${\Bbb C}$-vector spaces:

\[
J_{-1}(m):=Im[\Sigma Z_t H^0(I_X(m-1))] \subseteq %
J_0(m) \subseteq J_1(m) \subseteq \cdots \subseteq J_{n-1}(m) %
\subseteq J_n(m)=H^0(I_X(m)),
\]

\noindent
we chose a finite subset $ \{ F_{1,s,m},\ldots ,F_{k(s),s,m} \}$ from 
$J_s(m)$ which forms a ${\Bbb C}$-basis of $J_s(m)/J_{s-1}(m)$ and 
define closed subschemes $W_p \subset P$ 
and $W_p^* \subset P$ 
by the equations 
$\{ F_{1,s,m}, \cdots , F_{k(s),s,m}| 0 \leq s \leq p, \  m \in {\Bbb N}_0 \}$ 
and by 
$\{ F_{1,s,m}, \cdots , F_{k(s),s,m}| n \geq s \geq p, \ m \in {\Bbb N}_0  \}$, 
respectively for $p=0, \cdots ,n$.

\bigskip

%============Defn. of  Lefschetz Chain etc.==================

\begin{defn}[Lefschetz chain and dual Lefschetz chain]
Under the circumstances, we obtain two chains of closed 
subschemes $P$. The one is :

\[
j(X)=W_n \subseteq W_{n-1} \subseteq \cdots \subseteq %
W_0 \subseteq P 
\]
\noindent 
and is called {\em a Lefschetz chain} of $j(X)$. 
The other one is :

\[
j(X)=W_0^* \subseteq W_1^* \subseteq \cdots \subseteq %
W_n^* \subseteq P,
\]

\noindent
and is named {\em a dual Lefschetz chain} of $j(X)$. 
\end{defn}

\bigskip

Before claim our conjectures, we present 
fundamental properties of Lefschetz chains 
and dual Lefschetz chains.

%====== Thm. on Lefschetz chain etc. ==============

\begin{thm}
Let $X$ be a complex projective manifold of dimension 
$n>0$, 
$j: X \hookrightarrow P={\Bbb P}^{N}({\Bbb C})$ an 
arithmetically normal embedding.
% ,namely 
% $H^0(P,O_P(m)) \rightarrow H^0(X,O_X(m)$ is surjective 
% for all $m \in {\Bbb Z}$. 
Take a Lefschetz chain and 
a dual Lefschetz chain of $j(X)$ as above and fix them. 
The following properties holds.

\begin{conds}{1}

\condsitem  The submanifold $j(X)$ is a complete intersection 
if and only if the Lefschetz chain is of the form :

\[
j(X)=W_n \subset W_{n-1} = \cdots = W_0 = P. 
\]

\noindent
The similar equivalence holds on the dual Lefschetz chain by 
replacing the form:

\[
j(X)=W_0^* = W_1^* = \cdots = W_n^* \subset P.
\]

\condsitem Put $s=dim( Im[ L_X^n: \oplus_m H^0(N^{\vee}_X(m))%
\rightarrow \oplus_m H^n(N^{\vee}_X(m)\otimes \Omega^n_X) ] )$, 
then for the Lefschetz chain, the exact sequence:

\[
\begin{CD}
0 @>>> Im(N^{\vee}_{W_{n-1}}|_X) @>>> N^{\vee}_X @>>> %
N^{\vee}_{X/W_{n-1}} @>>> 0
\end{CD}
\]

\noindent
always splits and 
$N^{\vee}_{X/W_{n-1}} \cong \oplus^s O_X(-m_i)$.
Similarly for the dual Lefschetz chain, the exact 
sequence:

\[
0 @>>> N^{\vee}_{W_n^*}|_X @>>> N^{\vee}_X @>>> %
N^{\vee}_{X/W_n^*} @>>> 0
\]

\noindent
always splits and 
$N^{\vee}_{W_n^*}|_X \cong \oplus^s O_X(-m_i)$.

\condsitem Assume that the Standard Conjecture holds 
on the projective manifold $X$. For the Lefschetz chain,
if $W_p \not= W_{p-1}$, then there is an integer $m$ 
and a p-cycle $\xi$ such that 
$h^p \cdot \xi >0$ and 
$ \xi \cdot c_r(N^{\vee}_X(m)) \equiv_{num.eq.} 0$,
where $h=c_1(O_X(1))$. Also for the dual Lefschetz 
chain, if $W_p^* \not= W_{p+1}^*$, exactly the same holds.

\end{conds}
\end{thm}

\bigskip

\opf
For (2.5.1) and (2.5.2), we have only to apply Serre duality. 
The claim (2.5.3) is obtained by using the result of \cite{ABC} 
with a slight modification. To remove the condition "transverse 
to the zero section", we use Hironaka resolution for making the 
divisor normal crossing and study localized top Chern class instead 
of the zero locus of the section. 
\qed

\bigskip

Now we can describe our main conjectures as follows.

%============== Main Conjecture ====================

\begin{conj} \label{MCONJ}
Let $X$ be a complex projective manifold of dimension 
$n>0$, 
$j: X \hookrightarrow P={\Bbb P}^{N}({\Bbb C})$ an 
arithmetically normal embedding. Take a Lefschetz 
chain $\{ W_p \}_{p=0}^n$ and a dual Lefschetz 
chain $\{ W_p^* \}_{p=0}^n$ suitably. Then we expect 
the following properties hold by the suitable choice 
of the chains.

\begin{conds}{1}

\condsitem Each $W_p$ and $W_p^*$ are PG-shell of $j(X)$.

\condsitem Each $W_p$ and $W_p^*$ are reduced along $j(X)$.

\condsitem Each $W_p$ and $W_p^*$ are irreducible.

\condsitem Each restricted conormal sheaf 
$N^{\vee}_{W_p/W_{p-1}}|_X$ is a 
vector bundle on $X$ and is extendable to $W_{p-1}$ 
as a vector bundle. Similarly, each restricted 
conormal sheaf $N^{\vee}_{W_p^*/W_{p+1}^*}|_X$ 
is a vector bundle and is extendable to 
$W_{p+1}^*$ as a vector bundle. 

\condsitem Fix the manifold $X$ of dimension $n \geq 2$. 
Choose suitably the embedding $j$, 
Lefschetz chain $\{ W_p \}_{p=0}^n$, and a dual Lefschetz 
chain $\{ W_p^* \}_{p=0}^n$. Then a refinement of 
the Lefschetz chain or of the dual Lefschetz chain 
realizes the Working Hypothesis 
\ref{SWHP} {\em (cf. Theorem \ref{PBGRF})}.

\end{conds}

\end{conj}

\bigskip

%======= Prop. degree & meta-pent ==========

\begin{prop} \label{DEGMP}
Let $X$ be a connected complex projective manifold of dimension 
$n \geq 2$ and 
$j: X \hookrightarrow P={\Bbb P}^{N}({\Bbb C})$ 
an arithmetically normal embedding.
Assume that $j(X)$ is non-degenerate, namely 
no hyperplane in $P$ contains $j(X)$. 
For an equation 
$F \in H^0(P,I_X(m))$ of $j(X)$ 
which is an element of minimal generators of 
the ideal ${\Bbb I}_X$, take the class 
$[F] \in H^0(X,N^{\vee}_{X/P}(m))$ induced by the equation $F$.
Then we have:

\begin{conds}{1}
\condsitem If $m=2$, then the class 
$(\widehat{L}_X)^2([F]) \in H^{2}(X, %
N^{\vee}_{X/P}(2) \otimes \Omega_P^{2})$ is not zero.
More generally, take a non-zero class 
$\tau \in \widehat{gsyz}_X^q(q+1) \subseteq %
H^{1}(X, N^{\vee}_{X/P}(q+1) \otimes \Omega_P^{q})$ which naturally 
corresponds to an element of minimal generators of $q$-th syzygy 
module of $R_X$ in degree $q+1$, then 
$\widehat{L}_X(\tau) \in H^{2}(X, %
N^{\vee}_{X/P}(q+1) \otimes \Omega_P^{q+1})$ is not zero.

\condsitem If $m=3$ and the class 
$(\widehat{L}_X)^2([F]) \in H^{2}(X, %
N^{\vee}_{X/P}(3) \otimes \Omega_P^{2})$ is zero, 
then $q(X)=h^1(O_X) > 0$. 

\end{conds}
\end{prop}

\pf
Apply Theorem \ref{decomp} below.
\qed

\medskip

\begin{thm}[cf. \cite{MLO}] \label{decomp}
Let $X$ and $W$ be complex projective schemes of dimension 
$n \geq 0$ and of dimension $N$, respectively. 
Take an embedding 
$j: X \hookrightarrow W$ and assume that 
$j(X) \subset Reg(W)$.

For a non-zero integer $m$, 
the k-multiple of $\nu$-th meta-Lefschetz operator

\[
(\widehat{L}_{X/W}^{(\nu)})^k : H^{a}(X, %
I_{X}^{\nu + 1}/I_{X}^{\nu + 2}(m) \otimes \Omega_W^{b})
\rightarrow %
H^{a+k}(X, %
I_{X}^{\nu + 1}/I_{X}^{\nu + 2}(m) \otimes \Omega_W^{b+k})
\]

\noindent
is decomposed into 
$(\widehat{L}_{X/W}^{(\nu)})^k=c \cdot \overline{\delta}_{LFT}^{(\nu)} %
\circ (L_W)^{k-1} \circ \overline{d_I}$, where $c$ is 
a non-zero integer. 

\noindent
Similarly, the k-multiple of Lefschetz operator

\[
(L_W)^k : 
H^a(\Omega^{b}_{W}(m)|_{X_{(\nu)}})
\rightarrow %
H^{a+k}(\Omega^{b+k}_{W}(m)|_{X_{(\nu)}}),
\]

\noindent
is decomposed into 
$(L_W)^k=c' \cdot \overline{d_I}%
\circ (\widehat{L}_{X/W}^{(\nu)})^{k-1} \circ %
\overline{\delta}_{LFT}^{(\nu)}$, 
where $c'$ is a non-zero integer. 
\end{thm}

\medskip

%============ Remk meta-Lef. vs Lef.  =============

\begin{remk}
Proposition \ref{DEGMP} shows that the meta-Lefschetz operator 
has really finer information on the syzygies of the coordinate 
ring than the Lefschetz operator does. 
For example, take $X={\Bbb P}^{n}({\Bbb C})$ ($n \geq 2$), 
an embedding 
$j$=$\nu$-th Veronesean embedding ($\nu \geq 3$) and 
any equation $F$ of $j(X)$ in degree $2$, then it is 
easy to see that $({L}_X)^2([F]) \in H^{2}(X, %
N^{\vee}_{X/P}(2) \otimes \Omega_X^{2})$ is zero.
\end{remk}

We expect that the following theorem brings us 
a new idea for necessary refinements of Lefschetz 
chains and dual Lefschetz chains 
and helps us to solve our previous conjectures.

%====== Preservation of Filtration ==========

\begin{thm}[cf. \cite{MLO}] \label{PBGRF}
Let $X$ and $W$ be complex projective manifolds of dimension 
$n \geq 1$ and of dimension $N$, respectively. 
Take an embedding 
$j: X \hookrightarrow W$. 
Consider the exact sequence: 
$0 \rightarrow N^{\vee}_{X/W} \rightarrow %
\Omega_W^1|_X \rightarrow \Omega_X^1 \rightarrow 0$, 
which induces a canonical decreasing filtration $F_q^{\bullet}$: 
$\Omega_W^q|_X = F^0_q \supset F^1_q \supset %
\ldots \supset F^q_q \supset F^{q+1}_q= \{ 0 \} $ such that 
$F_q^p/F_q^{p+1} \cong \Omega_X^{q-p} \otimes %
\wedge^p  N^{\vee}_{X/W}$. 
A natural decreasing filtration on $H^{t,s}(m)=H^{s}(X, %
I_{X}^{\nu + 1}/I_{X}^{\nu + 2}(m)  \otimes \Omega_W^{t})$ 
is induced by putting 
$F^pH^{t,s}(m) = %
Im[H^s(X, I_{X}^{\nu + 1}/I_{X}^{\nu + 2}(m) %
\otimes F^p_t) \rightarrow  H^{t,s}(m)]$. 
Then the $\nu$-th meta-Lefschetz operator 
$\widehat{L}_{X/W}^{(\nu)}$ 
with respect to the embedding $j: X \hookrightarrow W$ 
keeps the filtration, namely 
$\widehat{L}_{X/W}^{(\nu)} (F^pH^{t,s}(m)) %
\subseteq F^pH^{t+1,s+1}(m)$.
\end{thm}

%%%%%%%%%%%%  Infinitesimal Methods. %%%%%%%%%%%%%%

\mysection{Infinitesimal Methods.}
\noindent

In this section, we introduce our simple 
tools which consist of two key concepts. 
These are mysteriously powerful for controlling higher 
obstructions appearing in 
the study of infinitesimal neighborhoods. 
These are important to consider the correspondence 
between subbundles of the normal bundle and intermediate 
ambient varieties.

%========= Defn of Diff. Split =================

\begin{defn}[Differential Splitting]

On a complex algebraic scheme $W$, 
we consider an ${\cal O}_W$-linear exact sequence of 
${\cal O}_W$-coherent sheaves:

\[
\begin{CD}
0 @>>>{ G} @>{\alpha}>> { F} @>{\beta}>> { E} @>>> 0.
\end{CD}
\]

\noindent
We say that this sequence {\em splits differentially of
order} $\leq {\mu}$ if there exists a (holomorphic 
${\Bbb C}$-linear) differential 
operator $\nabla_{\beta}: { E} \rightarrow { F}$ of order 
$\leq {\mu}$ such that $ \beta \circ \nabla_{\beta} = Id_{ E}$, 
namely, the operator $\nabla_{\beta}$ gives a splitting in the 
category of abelian sheaves. 
It is easy to see that
this condition is equivalent to the condition that the existence
of two differential operators 
$\nabla_{\alpha}:F \rightarrow G$ 
and $\nabla_{\beta}:E \rightarrow F$ 
of order $\leq {\mu}$ which satisfy : 
$ \beta \circ \nabla_{\beta} = Id_{ E}$ ; 
$ \nabla_{\alpha} \circ \alpha = Id_{ G}$ ; and 
$ \alpha \circ \nabla_{\alpha} + \nabla_{\beta} \circ %
\beta  = Id_{ F}$. When the scheme $W$ is smooth and 
the sheaf $E$ is of locally free, 
the condition of splitting 
differentially of some order 
is equivalent to the condition 
in terms of ${\cal D}_W$-modules 
that the sequence:

\[
\begin{CD}
0 @>>>{ G} \otimes {\cal D}_{W} @>{\alpha}>> { F} \otimes %
{\cal D}_{W} @>{\beta}>> { E} \otimes {\cal D}_{W} @>>> 0,
\end{CD}
\]

\noindent
splits in the category of right ${\cal D}_{W}$-modules, 
where ${\cal D}_W$ denotes the sheaf of holomorphic 
linear differential operators on $W$.
\end{defn}

As showed in \cite{ODJ}, 
there are many examples where differential splittings 
are observed. One of the typical examples is given as 
follows.

\bigskip

%========== Exampl. of Grassmann ===============

\begin{exam}
Let $V$ be a complex algebraic scheme, $E$ a vector bundle
on $V$, $f:G=Grass(E, r) \rightarrow V$ the Grassmann bundle which
parameterizes quotient r-bundles of $E$. Consider the universal sequence
on $G$:

\[
\begin{CD}
0 @>>> S @>{\alpha_G}>> f^{*}E @>{\beta_G}>> Q @>>> 0.
\end{CD}
\]

\noindent
Then this universal sequence splits differentially of order $ =1 $
(Obviously it never splits $O_G$-linearly). 
\end{exam}

\bigskip

%========= Defn. of H^p-GLC ====================

\begin{defn}[$H^{p}$-G.L.C.]
Let $W$ be a noetherian scheme, $X$ a closed 
subscheme of $W$ which is defined by a sheaf of ideals $I_X$,
$E$ a coherent ${\cal O}_W$-module. \\

\begin{conds}{1}

\condsitem For each non-negative integer $\mu$, we set the $\mu$-th 
infinitesimal neighborhood 
$X_{(\mu)}$ of $X$ in $W$ to be $(|X|, {\cal O}_{W}/I_{X}^{\mu+1})$ 
and the restricted sheaf $E_{(\mu)}$ of $E$ to $X_{(\mu)}$ to be 
$E/I_{X}^{\mu+1} E$ as usual. Let $\nu$ be a non-negative integer. 
We say that {\em the} $H^{p}$-{\em global lifting 
criterion of the coherent sheaf  E  holds at the (infinitesimal) 
lifting level} $\lambda$ 
{\em along} $(X_{(\nu)},X)$ if the equality:

\begin{eqnarray*}
&{}&Im[H^{p}(W,E) \rightarrow H^{p}(X_{(\nu)},E_{(\nu)})] \\
&{}&=Im[H^{p}(X_{(\nu + \lambda)},E_{(\nu + \lambda)}) %
\rightarrow  H^{p}(X_{(\nu)},E_{(\nu)})] 
\end{eqnarray*}

\noindent
holds in the space of $H^{p}(X_{(\nu)},E_{(\nu)})$. 
This condition is abbreviated as "$H^{p}$-G.L.C. of E holds at 
level $\lambda$ along $(X_{(\nu)},X)$".\\

\condsitem It is called that the $H^{p}$-global lifting criterion 
of the coherent sheaf  E  holds {\em uniformly} at the 
(infinitesimal) 
lifting level $\lambda$ along X if for any positive 
integer ${\nu}$, $H^{p}$-G.L.C. of E holds at 
level $\lambda$ along $(X_{(\nu)},X)$.
This condition is also abbreviated as "$H^{p}$-G.L.C. 
of E holds {\em uniformly} at lifting level $\lambda$ along $X$".
\end{conds}

\end{defn}

Let us show one of the results in \cite{PGLC} 
as the simplest example for showing the powerfulness 
of our previous two key concepts.

%=========== Thm on GLC of Q-type ============

\begin{thm}[Quotient Type] \label{QPGLC}
Let $W$ be a complex algebraic scheme.
For an exact sequence of $O_{W}$-coherent sheaves:

\[
0 \rightarrow G \stackrel{\alpha}{\rightarrow} F %
\stackrel{\beta}{\rightarrow} E \rightarrow 0 
\]

\noindent
connected by $O_{W}$-{\it linear} homomorphisms $\alpha$ and 
$\beta$, assume that this sequence splits differentially 
of order $\lambda$. 
If the $H^{p}$-lifting criterion on the sheaf $F$ holds
at the level $\mu$ along $(X_{(\nu)},X)$,
then the $H^{p}$-lifting criterion on the sheaf $E$ holds 
at the level $\lambda + \mu$ along $(X_{(\nu)},X)$. 

\end{thm}

\noindent
\pf It is enough to show that 
for any class 
$\phi \in H^p(X_{(\nu)}, E_{(\nu)})$ 
which is an image of a class of 
$H^p(X_{(\nu + \lambda)}, E_{(\nu + \lambda)})$,  
the class $\phi$ can be lifted to 
$H^p(W,E)$.

Let us consider six natural 
$O_{W}$-linear homomorphisms:  
$e : E \rightarrow E_{(\nu)}$, 
$\overline{e} : E_{(\lambda + \mu +\nu)} %
\rightarrow E_{(\nu)}$,\  $ r: E \rightarrow %
E_{(\lambda + \mu +\nu)}$,
$f : F \rightarrow F_{(\nu)}$, 
$\overline{f} : F_{(\mu + \nu)} %
\rightarrow F_{(\nu)}$,
\ and \ $ s: F \rightarrow F_{(\mu + \nu)} $,
which satisfy  $ e= \overline{e} \circ r $
 \ and \ $ f= \overline{f} \circ s $.
Since the differential operator $ \nabla : E \rightarrow F $
 is of ${\Bbb C}$-linear 
and of order $ \lambda $, it induces a homomorphism 
of abelian sheaves
$ \overline{\nabla} : E_{(\lambda + \mu + \nu)} %
\rightarrow F_{(\mu + \nu)}$ which 
satisfies $s \circ \nabla = \overline{\nabla} \circ r $. 
Then, using carefully the commutativities of the maps 
already checked, we see that:

\[
\overline{\beta} \circ \overline{f} \circ \overline{\nabla}%
 \circ r %
= \overline{\beta} \circ \overline{f} \circ s \circ \nabla %
= \overline{\beta} \circ f \circ \nabla = e \circ%
 \beta \circ \nabla %
= e \circ Id_{E}= \overline{e} \circ r,
\]

\noindent
where $ \overline{\beta} : F_{(\nu)} \rightarrow E_{(\nu)}$ 
denotes the natural $O_{W}$-linear homomorphism induced 
by $\beta : F \rightarrow E $.  
Considering all the homomorphisms given above as the 
homomorphisms 
in the category of abelian sheaves, the surjectivity of 
the homomorphism $r$ (at each stalk) implies that:

\[
\overline{\beta} \circ \overline{f} \circ %
\overline{\nabla} =\overline{e}.
\] 

\noindent
Now, by assumption, we can take a class  
$\overline{\psi}  \in H^{p}(X_{ ( \lambda + \mu + \nu) }, %
E_{(\lambda + \mu + \nu)})$ 
whose image by the map $\overline{e}$ 
coincides with the given class 
$\phi$ of $  H^{p}(X_{(\nu)}, E_{(\nu)}) $.
Then, taking $H^p$ of the sheaves introduced in the above, 
we have the following (a partially non-commutative) diagram:

\[
\begin{diagram}
\node{H^p(F)} 
	\arrow[3]{e,t}{\beta} 
	\arrow{s,l}{Id}
\node[3]{H^p(E)} 
	\arrow{s,r}{Id}
\\
\node{H^p(F)}
	\arrow{se,t}{s}
	\arrow{sse,b}{f}
\node[3]{H^p(E)}
	\arrow[3]{w,t}{\nabla}
	\arrow{sw,t}{r}
	\arrow{ssw,r}{e}
\\
\node[2]{H^p(F_{(\mu + \nu)})}
	\arrow{s,l}{\overline{f}}
\node{H^p(E_{(\lambda + \mu + \nu)})}
	\arrow{w,t}{\overline{\nabla}}
	\arrow{s,r}{\overline{e}}
\\
\node[2]{H^p(F_{(\nu)})}
	\arrow{e,b}{\overline{\beta}}
\node{H^p(E_{(\nu)}),}
\end{diagram}
\]

\noindent
By the assumption that 
$H^{p}$-G.L.C. of the sheaf F holds 
at the level $\mu$ along $(X_{(\nu)},X)$, the class 
$\overline{f} \circ \overline{\nabla}(\overline{\psi})$, 
which is the image of 
$\overline{\nabla}(\overline{\psi}) \in H^p(F_{(\mu + \nu)})$ 
by the map $\overline{f}$, 
can be lifted to a class  $\sigma \in H^{p}(W,F)$, 
namely $ f(\sigma)= \overline{f} \circ \overline{\nabla}( %
\overline{\psi} )$. 
 Then, putting
$\psi$ to be $ \beta ( \sigma ) $,  we see that:

\[ 
e(\psi) = e \circ \beta (\sigma) = %
\overline{\beta} \circ f ( \sigma) %
= \overline{\beta} ( \overline{f} \circ %
\overline{\nabla} ( \overline{\psi} ))%
= \overline{e} (\overline{\psi}) = \phi,
\]

\noindent
which is the desired conclusion.  \qed

\bigskip

%========= H^0-G.L.C. of O_P(m) =======================

\begin{cor}  \label{TLLF}
Let $V \subseteq P={\Bbb P}^N({\Bbb C})$ be 
a closed subscheme and 
$m \geq m_0$ non-negative integers. 
Assume that the restriction map 
$H^0(P,O_P(m_0)) \rightarrow H^0(V,O_V(m_0))$ 
is surjective. Then $H^0$-G.L.C. 
of $O_P(m)$ holds at level $m-m_0$ along $(V_{(0)},V)$. 
In other words, any section $\sigma \in H^0(V,O_V(m))$ 
can be lifted to $H^0(P,O_P(m))$ if and only if 
the section can be lifted to 
$H^0(V_{(m-m_0)}, O_{V_{(m-m_0)}}(m))$. 
\end{cor}

\smallskip

\pf 
By the assumption, $H^0$-G.L.C. of $O_P(m_0)$ holds at 
level $0$ along $(V_{(0)},V)$. 
We use induction on $m$ by starting from the case $m=m_0$. 
Take a positive integer $m > m_0$. 
We have only to apply Theorem \ref{QPGLC} 
to the Euler sequence:

\[
\begin{CD}
0 @>>> \Omega^1_P(m) @>>> \oplus O_P(m-1) @>>> O_P(m) @>>> 0,
\end{CD}
\]

\noindent
which splits differentially of order $1$ for positive 
integer $m$. (N.B. In case of $m=0$, this sequence never 
splits even in the sense of differential splitting.)
\qed

\bigskip

\begin{cor}  \label{TANG}
Let $V \subseteq P={\Bbb P}^N({\Bbb C})$ be 
a closed subvariety. Then $H^0$-G.L.C. 
of the tangent bundle $\Theta_{P}$ holds at 
level $2$ along $(V_{(0)},V)$. 
In other words, any section $\sigma \in H^0(V,\Theta_{P}|_{V})$ 
can be lifted to $H^0(P,\Theta_{P})$ if and only if 
the section can be lifted to 
$H^0(V_{(2)}, \Theta_{P}|_{V_{(2)}})$. 
\end{cor}
\pf 
By the assumption, $H^0$-G.L.C. of $O_P$ holds at 
level $0$ along $(V_{(0)},V)$.  Then Corollary \ref{TLLF} 
shows that $H^0$-G.L.C. of $O_P(1)$ holds at 
level $1$ along $(V_{(0)},V)$. 
Applying Theorem \ref{QPGLC} 
to the Euler sequence of the tangent bundle:

\[
\begin{CD}
0 @>>> O_P @>>> \oplus O_P(1) @>>> \Theta_P @>>> 0,
\end{CD}
\]

\noindent
which splits differentially of order $1$, we obtain the result.
\qed

%%%%%%%%%% Arithmetic Normality %%%%%%%%%%%%%%%%%%%

\mysection{Arithmetic Normality.}
\noindent

In this section, we discuss arithmetic normality 
from the two points of view. The first one is a 
viewpoint for clarifying our framework and strategy 
of studying the geometric structures of projective 
embeddings.
The second one is a viewpoint from 
Differential Geometry, 
which presents 
a criterion for arithmetic normality in terms 
of Differential Geometry. \\

For the first viewpoint, let us review weighted 
objects such as "weighted projections", which 
relates to "arithmetic normality" as a usual 
"projection" does to "linear normality".

%========== Weighted Projection =============

\begin{defn}[Weighted Projection]
For $N+L+1$-variables with weighted degree 
$wt.deg(Z_p)=s_p \geq 1$ $(p=0, \ldots ,N)$ ; 
$wt.deg(W_q)=w_q $ $(q=1, \ldots , L)$, take 
a weighted polynomial ring 
$T={\Bbb C}[Z_0, \ldots ,Z_N,W_1, \ldots ,W_L]$
and its polynomial subring $S={\Bbb C}[Z_0, \ldots ,Z_N]$.
By applying "Proj" operation, we get a rational map 
between the weighted projective spaces:

\[
\begin{diagram}
\node{ Proj(T)={\Bbb P}(s_0,\ldots ,s_N,w_1,\ldots ,w_L)} 
\arrow{e,t,..}{\pi} 
\node{ Proj(S)={\Bbb P}(s_0,\ldots ,s_N),}
\end{diagram}
\]

\noindent
which is called {\em a weighted projection} along the center 
$Z= \{ W_1=\ldots=W_L=0 \} $. 
\end{defn}

\begin{defn}[Weighted Linear Degeneration]
Consider a weighted polynomial ring 
$S={\Bbb C}[Z_0, \ldots ,Z_N]$ with $wt.deg(Z_p)=s_p$
and a closed subscheme 
$X \subset Proj(S)={\Bbb P}(s_0,\ldots ,s_N)=P$. 
If there is {\em a weighted linear} homogeneous polynomial $F \in S$, 
which is degree $1$ without weight in at least one variable,
e.g. $F=Z_0 + F_1(Z_1, \ldots Z_N)$, and if $X$ is a closed 
subscheme of the subscheme $Proj(S/(F)) \subset P$, then 
we say that the subscheme $X$ {\em degenerates weighted linearly}.
(In this case, the subscheme $X$ can be isomorphically projected 
through a suitable weighted projection.) 
\end{defn}

%======== Thm. Embed is w-proj & skew-arith-normal ======

\begin{lem} \label{FBWP}
Let $X$ be a complex projective scheme of 
dimension $n \geq 0$ and 
$j: X \hookrightarrow P={\Bbb P}^{N}({\Bbb C})={\Bbb P}(1,\ldots ,1)$ 
an embedding to a projective $N$-space (in a usual sense). 
Then there is a weighted projective space :
$Q={\Bbb P}(1^{N+1},w_1,\ldots ,w_L)$ and an embedding :
$\widetilde{\jmath}: X \hookrightarrow Q$ which make 
the commutative diagram:

\[
\begin{diagram}
\node{} \node{Q} \arrow{s,r,..}{\pi}\\
\node{X} \arrow{e,b}{j} \arrow{ne,t}{\widetilde{\jmath}} \node{P}
\end{diagram}
\] 

\noindent 
and satisfy the surjectivity on the natural map: 
$H^0(Q,O_Q(m)) \rightarrow H^0(X,O_X(m))$ 
for every non-negative integers $m$.
\end{lem}

\bigskip

Since several people asked me a proof for this lemma, 
it may be a little worth writing down its proof here.

\medskip

\pf
The idea is very simple and is only to add enough variables 
with suitable weighted degree. The argument goes as follows.
Let us put the vector space 
$V$ to be $Im[H^0(P,O_P(1)) \rightarrow H^0(X,O_X(1))]$ 
and the section $\sigma_t \in V$ to be 
the image of 
$Z_t \in H^0(P,O_P(1))$ for 
$t=0,1, \ldots , N$, where $N=dim(H^0(P,O_P(1)))-1$ and 
$ \{ Z_t \}_{t=0}^N $ form a ${\Bbb C}$-basis of $H^0(P,O_P(1))$. 
Since the line bundle $O_X(1)=j^*O_P(1)$ is ample, 
there are only finitely many positive integers \ $m$ \ such 
that 
$dim Coker[V \otimes H^0(X,O_X(m-1)) \rightarrow H^0(X,O_X(m))]%
=c_m \not= 0$.  
Set $ \{ m(1), \ldots ,m(u) \} = %
\{ m \in {\Bbb N} | c_m \not= 0 \}$ 
and $L=c_{m(1)} + c_{m(2)} + \cdots + c_{m(u)}$, where 
$1 \leq m(1) \leq \cdots \leq m(u)$.
Now we take the sections 
$\tau_1, \ldots , \tau_L$ such that 
$\tau_{ c_{m(1)}+ \cdots +c_{m(s-1)}+1 } ,\cdots , %
\tau_{ c_{m(1)}+ \cdots +c_{m(s)} } %
\in H^0(X,O_X(m(s)))$ induce the ${\Bbb C}$-basis of 
$Coker[V \otimes H^0(X,O_X(m(s)-1)) \rightarrow H^0(X,O_X(m(s)))]$ 
for $s=1, \ldots, u$. 
Take variables $W_k$ with $deg(W_k)=w_k$ corresponding to 
the section 
$\tau_k \in H^0(X,O_X(w_k))$ for $k=1, \ldots , L$, 
namely $w_k=m(s)$ if 
$ c_{m(1)}+ \cdots +c_{m(s-1)}+1 \leq k %
\leq  c_{m(1)}+ \cdots +c_{m(s)}$.
Now we have two essentially surjective ring 
homomorphisms :
$T={\Bbb C}[Z_0, \ldots, Z_N,W_1, \ldots, W_L] \rightarrow %
\widetilde{R_X}=\oplus_m H^0(X,O_X(m))$ and 
$S={\Bbb C}[Z_0, \ldots, Z_N] \rightarrow %
\widetilde{R_X}=\oplus_m H^0(X,O_X(m))$ 
by sending $Z_t$ to $\sigma_t$ 
and $W_k$ to $\tau_k$, which make a commutative diagram:

\[
\begin{diagram}
\node{} \node{T} \arrow{sw} \\
\node{\widetilde{R_X}} \node{S.} \arrow{w} \arrow{n,r}{inclusion}
\end{diagram}
\]

\noindent
Taking their "Proj", we obtain the result. 
(N.B. For simplicity, we constructed 
the ring $T$ rather roughly and it may have 
dispensable variables.)
\qed

%=================== Our Strategy =====================

\bigskip
Here we would like to make a discussion 
on a framework and a strategy for our research. 
Generally the weighted projective space $Q$ has singularities 
and the sheaf $O_Q(m)$ is not a line bundle but only 
a reflexive sheaf.
On the other hand, Lemma \ref{FBWP} above shows that any 
projective embedding is a composition of a weighted projection and 
an embedding into a weighted projective space 
which is very similar to an arithmetically 
normal embedding. 

Hence, to study 
the geometric structures of projective embedding, 
we can divide the problem into the three problem:
(a) investigate the arithmetically normal embeddings ;
(b) generalize the results of (a) into the case of 
weighted projective spaces (e.g. Working Hypothesis 
in weighted version);
(c) study the effects of weighted projections on the 
intermediate ambient varieties and on weighted G-shells
("weighted G-shell" is similarly defined by using 
$Tor^T_q(-,T/T_{+})$ instead of $Tor^S_q(-,S/S_{+})$). 

Relating to the problem (c) above, 
we should notice the fact that 
even if we have a good intermediate ambient variety 
$W$ with $\widetilde{\jmath}(X) \subset W \subset Q$, 
the variety $W$ may collapse by the weighted projection 
but the variety $X$ itself is projected isomorphically.
Thus we believe that 
the arithmetic normality is a natural condition as 
the fundamental assumption for our research in the 
first step, because we can ignore the difficulty 
arising from weighted projections. 

The arithmetic normality is equivalent 
to $H^0$-G.L.C. of $O_P(m)$ holding at level $0$ along 
$(X_{(0)},X)$ for every positive integer $m$ 
as we used it in the proof of Corollary \ref{TLLF}. 
Since the bundle $O_P(m)$ is a building block for coherent sheaves, 
the assumption of arithmetic normality makes the 
higher obstruction control much more easier than without it.
We might be going a bit too far, but 
the difficulty of higher obstruction can 
be sometimes explained 
by relating with weighted projections.

\bigskip

%=========== Review on Diff. Geom. ============
Now we proceed to the second viewpoint on arithmetic 
normality, namely that from Differential Geometry.
Let us recall the concepts of complex differential geometry.
Take a connected complex projective submanifold 
$X \subseteq P={\Bbb P}^{N}({\Bbb C})$ 
of dimension $n > 0$. By inducing 
a metric on $X$ from the Fubini-Study metric 
on $P$, we consider $X$ to be a K\"{a}hler manifold.
Consider the exact sequence of vector 
bundles with induced Hermitian metrics :

\[
\begin{CD}
0 @>>> N^{\vee}_{X/P} @>>> \Omega^1_P|_X @>>> \Omega^1_X @>>> 0.
\end{CD}
\]

\noindent
Then we have Hermitian connections 
$\nabla : {\cal A}^0(\Omega^1_P|_X) \rightarrow %
{\cal A}^1(\Omega^1_P|_X)$ and 
$\nabla_0 : {\cal A}^0(N^{\vee}_{X}) \rightarrow %
{\cal A}^1(N^{\vee}_{X})$, which induce a 
${\cal C}^{\infty}$-section 
$A=\nabla |_{N^{\vee}} - \nabla_0 \in %
{\cal A}^{(1,0)}(Hom(N^{\vee}_{X}, \Omega^1_X))$ 
of (1,0)-form with values in $Hom(N^{\vee}_{X}, \Omega^1_X)$.
% by using the natural 
% ${\cal C}^{\infty}$-isomorphism 
% $(N^{\vee}_{X})^{\bot} \cong \Omega^1_X$. 
Instead of $N^{\vee}_{X}$, considering 
$\Omega^1_X$ to be a ${\cal C}^{\infty}$-subbundle 
of $\Omega^1_P |_X$ by using the Hermitian metric, 
we have a ${\cal C}^{\infty}$-section 
$B \in {\cal A}^{(0,1)}(Hom(\Omega^1_X, N^{\vee}_{X}))$ 
of (0,1)-form with values in $Hom(\Omega^1_X, N^{\vee}_{X})$.

\bigskip

%======= Prop. on 2nd Fund. Form ==================

The following properties are well-known 
(cf. \cite{EPCA},\cite{LDG},\cite{DGCVB},\cite{SFF}).

\begin{prop}[Second Fundamental Forms]
Under the circumstances, 

\begin{conds}{1}
\condsitem $B$ is an adjoint of $-A$. In other words, 
for $\xi \in {\cal A}^0(N^{\vee}_{X})$ and 
$\eta \in {\cal A}^0(\Omega^1_X)$, the equality 
$h(A \xi,\eta) + h(\xi,B \eta)=0$ holds, where $h(-,-)$ denotes 
the Hermitian metric on $\Omega^1_P|_X$.

\condsitem Since $B$ is $\overline{\partial}$-closed, 
it defines a class $[B] \in H^1(X,\Theta_X \otimes N^{\vee}_{X})$, 
which coincides with the infinitesimal ring extension class of 

\[
\begin{CD}
0 @>>> N^{\vee}_{X} @>>> O_P/I_X^2 @>>> O_X @>>> 0.
\end{CD}
\]

\noindent
The class $\sigma_{II}(X)=[B]$ is called 
{\em the second fundamental 
form of type (0,1)} for $X$.

\condsitem $A \in H^0(Sym^2(\Omega^1_X) \otimes N^{\vee}_{X})$.
This class $A$ is called the {\em holomorphic second fundamental 
form} of $X$ and coincides with the differential of the Gauss map 
induced by the embedding.
Also a linear system is defined by considering it at general 
point of $X$ {\em (\cite{LDG},\cite{SFF})}.

\end{conds}
\end{prop}

\bigskip

Now we take a smooth irreducible divisor $D$ on $X$. Then we 
have an exact sequence : 

\[
\begin{CD}
0 @>>> \Theta_D @>>> \Theta_X |_D @>>> N_{D/X} @>>> 0,
\end{CD}
\]

\noindent
and a natural induced homomorphism : 
$r_D :H^1(X,\Theta_X \otimes N^{\vee}_{X}) \rightarrow %
H^1(N_{D/X}\otimes N^{\vee}_{X}|_D)$.

\bigskip

Using these notation,
we can describe 
a criterion for arithmetic normality, 
which was first obtained in \cite{PCHSU} 
by applying the view point of weighted 
projection. Here we explain an outline of 
another proof simplified by using 
the tools introduced in $\S 3$.

%======== Criterion on Arith Normal ============

\begin{thm}[Hoobler-Speiser-Usa]
Let $X \subseteq P={\Bbb P}^{N}({\Bbb C})$ 
be a connected complex projective submanifold 
of dimension $n \geq 2$. 
Assume that $q(X)=h^1(O_X)=0$.
Then the following two conditions are 
equivalent.

\begin{conds}{1}
\condsitem $X$ is arithmetically normal.

\condsitem For any positive integer $m$ and 
any generic smooth member 
$D \in |O_X(m)|$, $r_D(\sigma_{II}(X))=0$. 
\end{conds}

\end{thm}

\opf 
Showing arithmetic normality is the essential part. 
We apply induction on $m$. 
Take a section $\tau_D \in H^0(X,O_X(m))$ defining 
the divisor $D$. It is enough to see that the section 
$\tau_D$ lifts to $H^0(P,O_P(m))$. 
The assumption : $q(X)=0$ and $n \geq 2$ 
shows $H^1(N^{\vee}_X)=0$ through an easy application 
of Hodge Theory. 
By a direct computation on the exact sequence:

\[
\begin{CD}
0=H^1(N^{\vee}_X) @>{\times \tau_D }>> %
H^1(N^{\vee}_X(m)) @>{s_D}>> %
H^1(N^{\vee}_X|_D(m)), 
\end{CD}
\]

\noindent
we see that the class $r_D(\sigma_{II}(X)) \in %
H^1(N^{\vee}_X|_D(m))$ coincides with 
$s_D(\overline{\delta}^{(0)}_{LFT}(\tau_D))$. 
The assumption $r_D(\sigma_{II}(X))=0$ implies 
that the obstruction class: 
$\overline{\delta}^{(0)}_{LFT}(\tau_D) \in %
H^1(X,N^{\vee}_X(m))$ vanishes, which means that 
the section $\tau_D$ lifts to $H^0(X_{(1)},O_{X_{(1)}}(m))$.
Then we apply Corollary \ref{TLLF} to get the result. 
%For precise calculation on 
%$\overline{\delta}^{(0)}_{LFT}(\tau_D)$.
\qed

%%%%%%%%%%%%%%%%%%% References %%%%%%%%%%%%%%%%%%%%%%%%

\end{document}